\documentclass[11pt,a4paper,reqno]{amsart}
\usepackage[T1]{fontenc}
\usepackage{lmodern}
\usepackage{amssymb}
\usepackage[textwidth=150mm,textheight=236mm,centering]{geometry}
\usepackage{microtype}
\usepackage{amscd}
\usepackage[all]{xy}
\usepackage[hidelinks]{hyperref}
\hypersetup{
  pdftitle={The Riemann--Roch theorem as a consequence of duality},
  pdfauthor={Alberto Navarro Garmendia and José Navarro Garmendia}
}
\newcommand{\SH}{\mathbf{SH}}
\newcommand{\KGL}{\mathrm{KGL}}
\newcommand{\HQ}{\mathrm H_{\mathbb Q}}
\newcommand{\KQ}{\KGL_{\mathbf Q}}
\newcommand{\CH}{\mathrm{CH}}
\newcommand{\Q}{\mathbf Q}
\newcommand{\unit}{\mathbf 1}
\newcommand{\Th}{\operatorname{Th}}
\newcommand{\Td}{\operatorname{Td}}
\newcommand{\ch}{\operatorname{ch}}
\newcommand{\EE}{\mathbb{E}}
\newcommand{\FF}{\mathbb{F}}
\newcommand{\Hom}{\operatorname{Hom}}
\newcommand{\iHom}{\underline{\Hom}}
\newcommand{\Spec}{\operatorname{Spec}}
\newcommand{\Tseries}{\mathsf T}
\theoremstyle{plain}
\newtheorem{teo}{Theorem}[section]
\newtheorem{propo}[teo]{Proposition}
\newtheorem{coro}[teo]{Corollary}
\newtheorem{lema}[teo]{Lemma}
\theoremstyle{definition}
\newtheorem{defi}[teo]{Definition}
\theoremstyle{remark}
\newtheorem{obser}[teo]{Remark}

\title[Riemann--Roch and duality]{The Riemann--Roch theorem\\
  as a consequence of duality}
\author{Alberto Navarro Garmendia}
\address{Departamento de Matemática Aplicada, Universidad Politécnica de Madrid, Avenida Juan de Herrera 4, 28040 Madrid, Spain}
\email{alberto.navarro@upm.es}
\author{José Navarro Garmendia}
\address{Departamento de Matemáticas, Universidad de Extremadura, Avenida de Elvas s/n, 06006 Badajoz, Spain}
\email{navarrogarmendia@unex.es}
\date{}

\begin{document}
\begin{abstract}
We show how Grothendieck's Riemann--Roch theorem follows from duality
for a morphism between smooth proper schemes over a field.
Grothendieck duality allows us to identify the Spanier--Whitehead dual of a
smooth proper scheme and, for a cohomology theory with Chern classes, provides
a version of Poincaré duality. The Todd class appears when one compares the
two Poincaré isomorphisms. This gives a proof for the morphism itself,
without choosing a factorization.
\end{abstract}
\maketitle

\section{Introduction}

Grothendieck's Riemann--Roch theorem states the following.

\begin{teo}[Grothendieck--Riemann--Roch]\label{thm:grr}
Let $k$ be a field and let $f\colon X\to Y$ be a $k$-morphism between
smooth proper $k$-schemes. The square
\[
  \xymatrix@C=5pc@R=2.4pc{
    K_0(X) \ar[r]^{f_!} \ar[d]_{\Td(T_X)\,\ch}
      & K_0(Y) \ar[d]^{\Td(T_Y)\,\ch} \\
    \CH^*(X)_{\Q} \ar[r]_{f_*} & \CH^*(Y)_{\Q}
  }
\]
commutes. Equivalently, for every $\xi\in K_0(X)$,
\begin{equation*}
  \Td(T_Y)\,\ch(f_!\xi)
    = f_*\bigl(\Td(T_X)\,\ch(\xi)\bigr)
       \qquad\text{in }\CH^*(Y)_{\Q}.
\end{equation*}
\end{teo}
\noindent Here $T_X$ and $T_Y$ are the tangent bundles over $k$, $f_!$ is the
direct image on the Grothendieck group of perfect complexes induced by
$Rf_*$, and $f_*$ is the proper direct image of cycles.
For a line bundle $L$ on a smooth projective geometrically connected
curve $C$ of genus $g$, the theorem gives the classical formula
$\chi(C,L)=\deg L+1-g$.

In \emph{Récoltes et Semailles}, Grothendieck lists duality and the
``Yoga Riemann--Roch--Grothendieck'' as two of the principal themes of
his work~\cite[vol.~I, pp.~42, 44]{ReS}. The latter has indeed proved fruitful.
For example, it brought together $K$-theory and intersection theory through
the rational Chern character $K_0(X)_{\Q}\simeq\CH^*(X)_{\Q}$, a fact that
has motivated many constructions in the theory of motives (Bloch's higher
Chow groups relate higher algebraic $K$-theory to motivic cohomology
\cite{Bloch,RiouK,Soule}). The theorem itself has a topological counterpart
due to Atiyah and Hirzebruch~\cite{AH} and an arithmetic counterpart due to
Gillet and Soulé~\cite{GS}. The classical proof for a projective morphism
factors $f$ as a closed immersion $X\hookrightarrow\mathbf P^n_Y$ followed
by the projection $\mathbf P^n_Y\to Y$, and proves the formula for each of
these maps separately~\cite{BS}, \cite[\S15.2]{Fulton}. This strategy
has successfully applied  to other contexts~\cite[p.~246]{Gillet},
\cite[\S4.3]{GRS}.

The purpose of this note is to show how Riemann--Roch follows, almost
formally, from Grothendieck's other theme, duality.\footnote{Grothendieck
recalled that his work on Riemann--Roch had suggested to him that
global duality should be formulated for a proper morphism
\cite[p.~6]{Lafforgue}. Here we follow the connection in the other
direction.} Three forms of duality enter the argument. The first is Grothendieck duality
as part of the development of the formalism of the six operations. The
familiar dual of a finite-dimensional vector space has a categorical
formulation in the work of Dold and Puppe~\cite[\S\S1--3]{DoldPuppe},
computing such dual in stable homotopy is Spanier--Whitehead duality. This is the
second form of duality. The last is a version of Poincaré duality for a
cohomology theory with Chern classes (namely, an oriented cohomology
theory). Spanier--Whitehead duality and Poincaré duality follow, almost
formally, from Grothendieck duality.

Going back to Riemann--Roch, on the cohomology of the Spanier--Whitehead
duals, the usual direct images in $K$-theory and in Chow theory are thus
both given by precomposition with the same transpose $D(f)\colon DY\to DX$
\textup{(}see~\eqref{eq:push}\textup{)}. Since the Chern character
acts by postcomposition, a Riemann--Roch square on the duals
commutes by associativity alone \textup{(}Lemma~\ref{lem:transpose}\textup{)}. That is to say, there is a commutative Riemann-Roch square
\[
  \xymatrix@C=5pc@R=2.4pc{
    K_0(DX) \ar[r]^{f_!} \ar[d]_{\ch}
      & K_0(DY) \ar[d]^{\ch} \\
    \CH^*(DX)_{\Q} \ar[r]_{f_*} & \CH^*(DY)_{\Q}
  }
\]
with no need of Todd correction whose proof follows, almost tautologically, from Grothendieck's duality.
Passing through Poincaré duality, we recover the Riemann--Roch theorem.

This argument has been in the air for some time ad we claim no originality, apart from the exposition. Its topological antecedent is Dyer's proof by Spanier--Whitehead duality~\cite{Dyer}, and it is the
key to removing the projectivity hypothesis in the algebraic
case~\cite{NN}. We have long considered writing this note, but we refrained from doing so as it contains no new results and its setup is quite close to~\cite{NN}. However, given the current shift in mathematics—moving the emphasis from technical details and new proofs toward understanding and digesting the underlying principles—and the fact that these ideas form an essential part of~\cite{NN2}, we believe this note may now be of value to the community.

\section{Three forms of duality}

\subsection{Grothendieck duality and the exceptional functors}
In note~$46_2$ of \emph{Récoltes et Semailles}, Grothendieck presents the six operations as a
formalism for duality in cohomology, both global and local~\cite[pp.~2, 5]{Lafforgue}.

Stable motivic homotopy theory carries this formalism~\cite[Ex.~2.4.48 and Th.~2.4.50]{CD}.
For each scheme $X$ of finite type over $k$, $\SH(X)$ is closed symmetric
monoidal, with unit $\unit_X$. The first two operations are its tensor product
and internal Hom, written $\otimes$ and $\iHom_X$. Thus a map from $C\otimes A$
to $B$ is the same as a map from $C$ to $\iHom_X(A,B)$.

The other four operations relate the categories over two schemes. For a
morphism $g\colon X\to S$, the third is pullback $g^*\colon\SH(S)\to\SH(X)$
and the fourth is its right adjoint, direct image $g_*$. For smooth
$p\colon X\to S$, pullback also has a left adjoint $p_\sharp$. It sends the
suspension spectrum of a smooth $X$-scheme $U$ to that of $U$ viewed over $S$
by composition with $p$. In particular $p_\sharp\unit_X=X_+$, where
$X_+=\Sigma_T^\infty(X\sqcup *)$ in $\SH(S)$. For proper $g$, the direct
image $g_*$ in turn has a right adjoint, denoted by $g^!$.

The substance of Grothendieck duality here is the construction, for every
separated $g$ of finite type, of the last two operations, direct image with
proper support $g_!$ and its right adjoint $g^!$, the exceptional pullback.
They are functorial in $g$, and $g_!$ satisfies base change and the projection
formula.

For a vector bundle $V$ over $X$, denote by $\Th_X(V)$ the suspension spectrum of
$V/(V\setminus 0_X)$ over $X$, where $0_X$ is the zero section. Thom objects are
invertible in $\SH(X)$. We write $\Th_X(-V)$ for the inverse of $\Th_X(V)$.

\begin{teo}[Grothendieck duality]\label{thm:global}
Let $g\colon X\to S$ be a separated morphism of finite type between schemes of finite
type over $k$. The exceptional functors form an adjoint pair $g_!\dashv g^!$.
If $g$ is proper, $g_!\simeq g_*$. If $g$ is smooth with relative
tangent bundle $T_g$, the exceptional functors are given by
\begin{equation}\label{eq:purity}
 \begin{aligned}
    g_!M &\simeq g_\sharp\bigl(\Th_X(-T_g)\otimes M\bigr),\qquad &&M\in\SH(X),\\
    g^!B &\simeq \Th_X(T_g)\otimes g^*B, &&B\in\SH(S).
 \end{aligned}
\end{equation}
\end{teo}
For the smooth formulas see Ayoub~\cite[\S\S1.5.3.1--1.5.3.2]{Ayoub}. The extension
from the quasi-projective setting to separated morphisms of finite type is given in
\cite[Ex.~2.4.48, Th.~2.4.50(2), (3) and Rem.~2.4.51]{CD}.

For an open immersion $j\colon U\hookrightarrow S$, $j^*$ is restriction and
$j_\sharp\simeq j_!$ is extension by zero. Thus $j_!M$ restricts to $M$ on $U$ and
vanishes on its closed complement~\cite[\S2.3.1 and Th.~2.4.50(3)]{CD}.

\begin{obser}[Global and local duality]\label{rem:lafforgue}
We follow Grothendieck's terminology, as recalled by Lafforgue~\cite[pp.~2--3]{Lafforgue}.
\emph{Global duality} expresses the exceptional adjunction as an isomorphism in $\SH(S)$,
\begin{equation}\label{eq:global}
 \iHom_S(g_!M,B)\simeq g_*\iHom_X(M,g^!B).
\end{equation}
Indeed, a map $C\otimes g_!M\to B$ corresponds, by the projection formula and
$g_!\dashv g^!$, to a map $g^*C\otimes M\to g^!B$. By tensor--Hom adjunction
and $g^*\dashv g_*$, this is a map from $C$ to the right-hand
side~\cite[Prop.~2.3.53]{Ayoub}, \cite[Th.~2.4.50(5)]{CD}.

Put $D_KM=\iHom(M,K)$. Following Riou's definition for constructible
complexes~\cite[Def.~0.1]{Dualite}, we call an object $K$ of a full subcategory
$\mathcal C$ \emph{dualizing} if $D_K$ preserves $\mathcal C$ and the canonical
map $M\to D_K(D_KM)$ is an isomorphism for every $M\in\mathcal C$.
Grothendieck calls this property \emph{local~duality}, or \emph{biduality}.

For smooth $p\colon X\to\Spec k$, the object $\omega_p=p^!\unit_k\simeq\Th_X(T_X)$
is invertible, so $D_{\omega_p}$ interchanges $\unit_X$ and $\omega_p$. Since
$p_!\omega_p\simeq p_\sharp\unit_X=X_+$, taking $M=\omega_p$ and $B=\unit_k$
in~\eqref{eq:global} gives $\iHom_k(X_+,\unit_k)\simeq p_*\unit_X$.
\end{obser}

\goodbreak
\subsection{The Spanier--Whitehead dual}
\begin{defi}[Spanier--Whitehead dual]\label{def:sw}
The \emph{dual} of an object $A$ of $\SH(S)$ is
\[
                       DA=\iHom_S(A,\unit_S).
\]
The \emph{transpose} of a map $h\colon A'\to A$ is the map
$Dh\colon DA\to DA'$ induced by precomposition with $h$. The evaluation map $DA\otimes A\to\unit_S$ induces a map
\[
        DA\otimes B\longrightarrow\iHom_S(A,B).
\]
One says that $A$ is \emph{strongly dualizable} if this map is an isomorphism for every $B$. In that case evaluation also identifies $A$ with its double dual.
\end{defi}

For a smooth $k$-scheme $X$, write $DX=D(X_+)$. Riou proves the following result for smooth projective schemes~\cite[Th.~2.2]{RiouSW}.\footnote{The
description of the dual in terms of the negative tangent bundle is also called
\emph{Atiyah duality}, after its topological counterpart~\cite[Th.~3.3]{Atiyah}.}

\begin{propo}[Motivic Spanier--Whitehead duality]\label{prop:sw}
Let $p\colon X\to\Spec k$ be smooth and proper. There is a canonical isomorphism
in $\SH(k)$
\begin{equation*}
                          DX\simeq p_\sharp\Th_X(-T_X).
\end{equation*}
\end{propo}
\begin{proof}
Let $C\in\SH(k)$. Since $C\otimes X_+\simeq p_\sharp p^*C$
by the projection formula, the adjunctions $p_\sharp\dashv p^*\dashv p_*$ give
\[
\begin{aligned}
 \Hom_{\SH(k)}(C,DX)
    &\simeq \Hom_{\SH(k)}(C\otimes X_+,\unit_k)\\
    &\simeq \Hom_{\SH(X)}(p^*C,\unit_X)\\
    &\simeq \Hom_{\SH(k)}(C,p_*\unit_X).
\end{aligned}
\]
Thus $DX\simeq p_*\unit_X$ by Yoneda. Properness gives
$p_*\unit_X\simeq p_!\unit_X$, and~\eqref{eq:purity} gives
$p_!\unit_X\simeq p_\sharp\Th_X(-T_X)$. For smooth proper $X$, the spectrum
$X_+$ is moreover strongly dualizable~\cite[\S3, p.~3616]{Hoyois}, but we will not use this fact.
\end{proof}

\subsection{Poincaré duality for an oriented theory}
A commutative motivic ring spectrum $\EE$ plays the role of a \emph{ring of coefficients of cohomology theory}
(cf.~\cite[p.~1]{Lafforgue}). The cohomology of a spectrum $A$ with these
coefficients is
\[
             \EE^{a,b}(A)=\Hom_{\SH(k)}\bigl(A,\EE(b)[a]\bigr).
\]
For a smooth $k$-scheme $X$ we write $\EE^{a,b}(X)=\EE^{a,b}(X_+)$ and $\EE(X)=\bigoplus_{i,j}\EE^{i,j}(X)$.
The six operations give the basic properties of cohomology, including pullbacks,
cup products~\cite[Rem.~1.5]{NN}.

Suppose now that $\EE$ is oriented, that is to say, we fix the first Chern class of the universal line bundle on $\mathbf{P}^\infty$; hence, all Chern classes for every vector bundle. Let $V$ be a vector bundle of rank $r$ on
$X$ and let $q\colon \overline V=\mathbf P(V\oplus\mathcal O_X)\to X$ be its projective
completion.  Write $s\colon X\to V$ for the zero section.
We write $\Th_X(V):=V/(V-s(X))$. As in topology, excision and homotopy provides and $\mathbf A^1$-equivalence, 
$\overline V/\mathbf P(V)\ \simeq \ \Th_X(V)$ in $\SH(X)$~\cite[\S3, Prop.~2.17(3)]{MV}. For $p\colon X\to\Spec k$, the the projective bundle formula give a split exact sequence
of bigraded $\EE(X)$-modules
\[
 0\longrightarrow\EE(\Th_X(V))\longrightarrow
 \EE(\overline {V})\longrightarrow\EE(\mathbf P(V))\longrightarrow0 .
\]
Hence, $\EE(\Th_X(V))$ is a free $\EE(X)$-module of rank 1 generated by a class $t^{\EE}(V)$ named \emph{Thom class}, of bidegree $(2r,r)$. It is determined by the fact that it maps to the last Chern class of the universal quotient bundle of $\overline{V}$. Multiplication by $t^{\EE}(V)$ is the Thom isomorphism. Thom classes multiply
in short exact sequences and extend to virtual bundles by taking inverse
Thom classes~\cite{RiouK}. The Thom isomorphism for $-T_X$, together with
Proposition~\ref{prop:sw} gives Poincaré duality.

\begin{coro}[Poincaré duality]\label{cor:pd}
For a smooth proper $k$-scheme $X$ of pure dimension $d$ and an oriented theory $\EE$,
there is a canonical isomorphism
\begin{equation}\label{eq:pd}
 \xymatrix@C=4pc@R=0.4pc{
    \EE^{a,b}(X) \ar[r]^-{\sim} & \EE^{a-2d,b-d}(DX) \\
    \alpha \ar@{|->}[r] & \alpha\cdot t(-T_X).}
\end{equation}
\end{coro}

\section{The direct image by transposition}
Let $f\colon X\to Y$ be a morphism between smooth proper $k$-schemes.
Its transpose $D(f)\colon DY\to DX$ makes the cohomology of the duals
covariant: a class $a\colon DX\to\EE(s)[r]$ is sent to
$a\circ D(f)\colon DY\to\EE(s)[r]$.
For an oriented theory the commutative square
\begin{equation}\label{eq:push}
 \vcenter{\hbox{$\xymatrix@C=5pc@R=2.2pc{
    \EE^{p,q}(X) \ar@{-->}[r]^-{f_*} \ar[d]_-{\text{Poincaré}}^-{\sim}
       & \EE^{p+2c,q+c}(Y) \ar[d]^-{\text{Poincaré}}_-{\sim} \\
    \EE^{r,s}(DX) \ar[r]_-{D(f)^*} & \EE^{r,s}(DY) }$}}
\end{equation}
defines the direct image $f_*$. Here $c=d_Y-d_X$,
$r=p-2d_X$ and $s=q-d_X$. Notice that the lower arrow is precomposition with the same
$D(f)$ for every $\EE$, and that we obtain it from $f$ itself, without the need of factoring $f$.

\begin{obser}
The lower arrow in~\eqref{eq:push} is proper pushforward in Borel--Moore
homology, so the upper arrow is the direct image of~\cite[Def.~2.2]{NN}.
In particular, on $K_0$ it is induced by $Rf_*$, and on rational Chow groups it is the usual
proper pushforward of cycles~\cite[Cor.~2.12]{NN}, \cite[Ex.~1.2.10(1)]{DegBiv}.
\end{obser}

\begin{obser}\label{rem:transfer}
Let us rephrase the above construction. Put $T_f=T_X-f^*T_Y$.
Since $f_*=f_!$, the unit $B\to f_*f^*B$ induces
$p_{Y!}B\to p_{X!}f^*B$. Taking $B=\Th_Y(T_Y)$ and applying smooth purity
gives the transfer
\[
     \tau_f\colon Y_+\longrightarrow p_{X\sharp}\Th_X(-T_f).
\]
This map encodes the direct image in cohomology. For $\alpha\in\EE^{p,q}(X)$
and $c=\operatorname{rk}(-T_f)$,
\[
       f_*(\alpha)=\tau_f^*\bigl(\alpha\,t^{\EE}(-T_f)\bigr)
                       \quad\text{in }\EE^{p+2c,q+c}(Y).
\]
\end{obser}

\section{The Riemann--Roch theorem}

First of all, let us remark a simple observation: the maps defined in the previous section at the level of duals commute with  arrows between cohomology theories. To be precise:

\begin{lema}\label{lem:transpose}
Let $\phi\colon\EE\to\FF$ be a map in $\SH(k)$ and let $f\colon X\to Y$ be a
morphism of smooth $k$-schemes. In every bidegree $(r,s)$ there is a commutative square
\[
\xymatrix@C=4pc@R=1.6pc{
 \EE^{r,s}(DX) \ar[r]^{D(f)^*} \ar[d]_{\phi}
    & \EE^{r,s}(DY) \ar[d]^{\phi} \\
 \FF^{r,s}(DX) \ar[r]_{D(f)^*} & \FF^{r,s}(DY).}
\]
\end{lema}
\begin{proof}
For $a\colon DX\to\EE(s)[r]$, both paths give the composite
\[
 DY\xrightarrow{\ D(f)\ }DX\xrightarrow{\ a\ }\EE(s)[r]
                       \xrightarrow{\ \phi(s)[r]\ }\FF(s)[r]. \qedhere
\]
\end{proof}

To pass from the duals to the schemes themselves, it remains to see what
$\phi$ does to Poincaré's duality. From now on $\EE$ and $\FF$ are two coefficient rings  with chosen orientations, and $\phi$ is a unital map of ring
spectra. For the tautological line $\mathcal O(-1)$ on $\mathbf P^\infty$, put
$x=c_1^{\FF}(\mathcal O(-1))$. The projective bundle formula and the normalization of both
orientations on $\mathbf P^1$ determine a series
$\Tseries(t)=1+\sum_{n\geq1}a_n t^n$, with $a_n\in\FF^{-2n,-n}(k)$, by
\[
        \phi\bigl(c_1^{\EE}(\mathcal O(-1))\bigr)=\frac{x}{\Tseries(x)}.
\]
Write $\Tseries_\times$ for its multiplicative extension to virtual bundles.
Thus $\Tseries_\times(L)=\Tseries(c_1^{\FF}(L))$ for every line bundle $L$
and $\Tseries_\times(-\xi)=\Tseries_\times(\xi)^{-1}$.

\begin{lema}\label{lem:virtual-thom}
For a virtual vector bundle $\xi$ on a smooth $k$-scheme $X$,
\[
       \phi\bigl(t^{\EE}(\xi)\bigr)
              =\Tseries_\times(-\xi)\,t^{\FF}(\xi).
\]
\end{lema}
\begin{proof}
Write $\xi=V-W$. Both sides are multiplicative in $\xi$, and $t(-W)$ is
the inverse Thom class of $W$. It suffices, therefore, to consider a vector
bundle $V$. By the splitting principle, we reduce the assertion to a line bundle $L$.

By homotopy invariance and the classification of line bundles by
$\mathbf P^\infty\simeq B\mathbf G_m$, it remains to treat the tautological
line~\cite[\S4, Props.~3.7--3.8]{MV}, and this follows from the definition of $\Tseries$. 
\end{proof}

\begin{coro}\label{cor:orient}
For a smooth proper $k$-scheme $X$ of pure dimension $d$, one has the commutative square
\begin{equation}\label{eq:orient}
\vcenter{\hbox{$\xymatrix@C=5pc@R=1.9pc{
 \EE^{p,q}(X) \ar[r]^-{\text{Poincaré}}_-{\sim}
                 \ar[d]_{\Tseries_\times(T_X)\,\phi}
       & \EE^{p-2d,q-d}(DX) \ar[d]^{\phi} \\
 \FF^{p,q}(X) \ar[r]^-{\text{Poincaré}}_-{\sim}
       & \FF^{p-2d,q-d}(DX). }$}}
\end{equation}
\end{coro}
\begin{proof}
Poincaré duality sends $\alpha$ to $\alpha\,t(-T_X)$, so that applying
Lemma~\ref{lem:virtual-thom}  gives
\[
 \phi\bigl(\alpha\,t^{\EE}(-T_X)\bigr)
     =\phi(\alpha)\,\Tseries_\times(T_X)\,t^{\FF}(-T_X). \qedhere
\]
\end{proof}

\goodbreak
\begin{teo}\label{thm:coefficients}
Let $\phi\colon\EE\to\FF$ be a homomorphism of coefficient rings in $\SH(k)$,
each equipped with an orientation, and let $\Tseries$ be its associated series.
For a morphism $f\colon X\to Y$ between smooth proper $k$-schemes of pure
dimensions $d_X,d_Y$, put $c=d_Y-d_X$. Then one has the commutative square
\[
\xymatrix@C=4.8pc@R=1.8pc{
 \EE^{p,q}(X) \ar[r]^-{f_*^{\EE}} \ar[d]_{\Tseries_\times(T_X)\,\phi}
       & \EE^{p+2c,q+c}(Y) \ar[d]^{\Tseries_\times(T_Y)\,\phi} \\
 \FF^{p,q}(X) \ar[r]_-{f_*^{\FF}} & \FF^{p+2c,q+c}(Y).}
\]
Equivalently, for $T_f=T_X-f^*T_Y$ and $\alpha\in\EE^{p,q}(X)$,
\begin{equation}\label{eq:relative-rr}
  \phi(f_*^{\EE}\alpha)
       =f_*^{\FF}\bigl(\Tseries_\times(T_f)\,\phi(\alpha)\bigr).
\end{equation}
\end{teo}
\begin{proof}
Expanding the two direct images by~\eqref{eq:push} gives
\[
\xymatrix@C=4pc@R=1.9pc{
 \EE(X) \ar[r]^-{\text{Poincaré}}_-{\sim} \ar[d]_-{\Tseries_\times(T_X)\,\phi}
 & \EE(DX) \ar[r]^-{D(f)^*} \ar[d]_-{\phi}
 & \EE(DY) \ar[d]^-{\phi}
 & \EE(Y) \ar[l]_-{\text{Poincaré}}^-{\sim} \ar[d]^-{\Tseries_\times(T_Y)\,\phi} \\
 \FF(X) \ar[r]^-{\text{Poincaré}}_-{\sim}
 & \FF(DX) \ar[r]_-{D(f)^*}
 & \FF(DY) & \FF(Y) \ar[l]_-{\text{Poincaré}}^-{\sim} .}
\]
The outside squares commute by Corollary~\ref{cor:orient} and the middle square commutes
by associativity, $\phi\circ(a\circ D(f))=(\phi\circ a)\circ D(f)$.
\end{proof}

\subsection{Two classical formulas}
For the Chern character take $K=\KQ$ and let
$H=\bigoplus_{m\in\mathbf Z}\HQ(m)[2m]$ denote periodic Beilinson motivic cohomology. The motivic Chern character $\ch\colon K\xrightarrow{\sim}H$
induces the usual map $K_0(X)_{\Q}\to\CH^*(X)_{\Q}$ on a smooth $k$-scheme $X$ in degree $(0,0)$%
~\cite[Def.~6.2.3.9]{RiouK}. For the universal line $L$, writing $x=c_1^{\CH}(L)$ gives
\[
      \ch(1-[L^\vee])=1-e^{-x},\qquad
      \Td(x)=\frac{x}{1-e^{-x}}.
\]
Thus the series associated to $\ch$ in these coordinates is the Todd series. We deduce from
Theorem~\ref{thm:coefficients} the classical formula.

\begin{coro}[Grothendieck--Riemann--Roch]\label{cor:classical}
For a morphism $f\colon X\to Y$ between smooth proper $k$-schemes and
$\xi\in K_0(X)$,
\[
      \ch(f_!\xi)=f_*\bigl(\Td(T_f)\,\ch(\xi)\bigr)
                      \quad\text{in }\CH^*(Y)_{\Q}.
\]
\end{coro}

For $j\in\mathbf Z\setminus\{0\}$, the Adams operation lifts to a unital map of
ring spectra $\Psi^j\colon\KGL[1/j]\to\KGL[1/j]$~\cite[Def.~5.3.2 and the paragraph following it]{RiouK}.
Indeed, tensor powers preserve products and send the Bott element $\beta$ to $j\beta$.

\goodbreak
\begin{coro}[Adams--Riemann--Roch]\label{cor:adams}
For $f$ and $\xi$ as above and $j\ne0$,
\[
       \psi^j(f_!\xi)=f_!\bigl(\theta^j(T_f^\vee)^{-1}\,\psi^j(\xi)\bigr)
                         \quad\text{in }K_0(Y)\otimes\mathbf Z[1/j].
\]
\end{coro}
\begin{proof}
In degree-zero coordinates the series associated to $\psi^j$ is
$\frac{t}{1-(1-t)^j}$. If we write $\Theta_j$ for its multiplicative extension,
the statement follows from Theorem~\ref{thm:coefficients} if we prove that
\[
                    \Theta_j(\xi)=\theta^j(\xi^*)^{-1}.
\]

It is enough to check this equality for a line bundle $L\to X$.
For $j>0$, writing $L^*$ for the class of the dual line in $K_0(X)$, we have
\[
 \Theta_j(L)=\frac{c_1^K(L)}{1-(1-c_1^K(L))^j}
 =\frac{1-L^*}{1-(L^*)^j}
 =\frac{1}{1+L^*+\cdots+(L^*)^{j-1}}
 =\theta^j(L^*)^{-1}.
\]
For $j<0$ the same calculation uses $\theta^j(L^*)=B_j(L^*)$, where
$B_j(z)=(1-z^j)/(1-z)\in\mathbf Z[z,z^{-1}]$ and $B_j(1)=j$. 
\end{proof}

\end{document}